\documentclass[12pt]{amsart}
\usepackage{amsmath} 
\usepackage{amsfonts}
\newtheorem{theorem}[subsection]{Theorem}
\newtheorem{proposition}[subsection]{Proposition}
\newtheorem{corollary}[subsection]{Corollary}
\newtheorem{lemma}[subsection]{Lemma}

\theoremstyle{definition}
\newtheorem{definition}[subsection]{Definition}

\newtheorem{notation}[subsection]{Notation}

\theoremstyle{remark}
\newtheorem{remark}[subsection]{Remark}

\newtheorem{remarks}[subsection]{Remarks}

\def\la{\langle}
\def\ra{\rangle}

\def\cA{\mathcal{A}}
\def\cF{\mathcal{F}}

\def\cB{\mathcal{B}}

\def\ct{\mathbf{t}}
\def\cI{\mathcal{I}}
\def\cV{\mathcal{O}}
\def\HH{\mathcal{H}}
\def\CC{\mathbb{C}}
\def\NN{\mathbb{N}}
\def\RR{\mathbb{R}}
\def\QQ{\mathbb{Q}}

\def\ff{\varphi}

\def\hN{{(N)}}
\def\vN{\text{vN}}
\def\ii{{\bf i}}
\def\jj{{\bf j}}
\begin{document}
\sloppy

\title{On the structure of non-commutative white noises}
\author[C. K\"ostler]{Claus K\"ostler}
\author[R. Speicher]{Roland Speicher $^{(\dagger)}$}
\thanks{$^\dagger$ 
Research supported by a Discovery Grant and a Leadership Support
Initiative Award from the Natural Sciences and Engineering 
Research Council of Canada and by a Premier's Research Excellence
Award from the Province of Ontario}
\address{Queen's University, Department of Mathematics and Statistics,
Jeffery Hall, Kingston, ON K7L 3N6, Canada}

\email{koestler@mast.queensu.ca, speicher@mast.queensu.ca}

\begin{abstract}
We consider the concepts of continuous Bernoulli systems and 
non-commutative white noises. We address the question of 
isomorphism of continuous Bernoulli systems and show that for large 
classes of quantum L{\'e}vy processes one can make quite 
precise statements about the time behaviour of their moments. 
\end{abstract}
 
\maketitle

\section{Introduction}
 
In recent years it has become evident that looking on non-commutative
algebras (in particular, operator algebras) from a stochastic  point
of view can be quite fruitful. So the impressive progress  on our
understanding of the free group von Neumann algebras relies  on
Voiculescu's free probability approach \cite{Voi,VDN} and  the work of
Pisier and Xu on non-commutative martingales has  opened a new
direction of research \cite{PX}. In particular,  it seems that
non-commutative versions of classical stochastic  processes yield
interesting examples of non-trivial operator-algebraic  structures. At
the moment we are only at the beginning of an   understanding of the
richness of the world of non-commutative  processes; the present paper
is a contribution to a systematic  theory of such non-commutative
processes.
 
L\'evy processes, i.e.\ processes with stationary and independent
processes, or `white noises' as models for their `derivatives',  form
one of the most important  classes of classical stochastic processes
and the understanding  of their structure was instrumental for many
developments in  classical probability theory.  It is to be expected
that the understanding of non-commutative  versions of L\'evy
processes will be  an important step towards a deeper theory of
non-commutative  stochastic processes.
 
An axiomatic frame for the treatment of non-commutative white noises
was started by K\"ummerer \cite{Kum1,Kum2} and is further elaborated
by one of us \cite{Kos,Kos1,HKK}. Here we will address some of the
canonical basic   questions of this theory: namely how we can
distinguish  between  different non-commuative white noises; and what
can be said about the time behaviour of their moments. Even though a
general answer to these  problems for the class of all non-commutative
L\'evy processes seems to be out of reach (and might not even exist in
this generality, see  Section \ref{subsection:example-non-order}), we
are able to  provide  answers to these questions for some quite large
classes of non-commutative white noises.
  
We also want to point out that an It\^o integration theory  for
non-commutative L\'evy processes was established in \cite{HKK}.
However, our results here will not rely on this integration  theory.

\section{Basic Definitions}
We want to generalize the notion of a classical process with
stationary and independent increments to a non-commutative setting. In
the classical setting, it is not only the process itself which is of
importance, but sometimes one is more interested in the structure of
the associated filtration  of $\sigma$-algebras of the increments. In
the same way, we find it advantageous in the non-commutative case to
distinguish between the filtration generated by the process, and the
process itself. In the non-commutative setting,  the filtration is
given by the von Neumann algebras generated by (or, in the case of
unbounded operators, affiliated to) the increments of the processes --
we will address this data as a continuous Bernoulli system.  We will
restrict here to the finite case, i.e. where the underlying state
$\ff$ is a trace. More general definitions are possible (and desirable
for a general theory), however, here we want to consider only the
simplest non-trivial case.
 
\begin{notation} 
By $\cI$ we denote the set of all intervals $I\subset\RR$  of the form
$I=[s,t)$ for $-\infty< s<t<\infty$.  For $I=[s,t)\in\cI$ and
$u\in\RR$ we denote by $I+u$ the  interval
$$I+u:=[s+u,t+u).$$ 
\end{notation} 

\begin{definition}
A \emph{(non-commutative) continuous Bernoulli system}
$(\cA,\ff,(\cA_I)_{I\in\cI})$ consists of
\begin{enumerate}
\item[(i)]  a non-commutative probability space $(\cA,\ff)$, where
$\cA$ is  a von Neumann algebra with separable predual and $\ff$ is  a
faithful and normal trace $\ff$ on $\cA$;
\item[(ii)]  a filtration $(\cA_I)_{I\in\cI}$, where $\cA_I$ are von
Neumann subalgebras of $\cA$ such that the following  properties are
satisfied:
\begin{enumerate} 
\item 
\emph{global minimality}:
$$\cA=\bigvee\{\cA_I\mid I\in\cI\};$$
\item 
\emph{isotony}:
$$\cA_I\subset \cA_J\qquad\text{whenever $I\subset J$};$$
\item 
\emph{$\CC$-independence}: for all $I,J\in\cI$ with $I\cap
J=\emptyset$ we have that
$$\ff(ab)=\ff(a)\ff(b)$$  for all $a\in\cA_I$ and all $b\in\cA_J$.
\end{enumerate} 
\end{enumerate} 
\end{definition}

\begin{remarks} 
1) Notice that we can phrase the $\CC$-independence also in the
following form:  for all $s<t<u$,
$$\begin{matrix} 
\cA_{[s,t)}&\subset&\cA_{[s,u)}\\ 
\cup &\quad& \cup\\ 
\CC&\subset&\cA_{[t,u)} 
\end{matrix}$$ 
is a (not necessarily non-degenerate) commuting square \cite{Popa}. 
If the von Neumann algebra $\cA$ is commutative, $\CC$-independence
is equivalent to the usual notion of stochastic independence in
probability theory.  

2) In our tracial frame, $\CC$-independence is clearly equivalent  to
pyramidal independence:  for all $I,J\in\cI$ with $I\cap J=\emptyset$
we have that
$$\ff(a_1ba_2)=\ff(a_1a_2)\ff(b)\qquad\text{for all $a_1,a_2\in \cA_I$
and all  $b\in\cA_J$.}$$   In a more general, non-tracial frame, one
needs the existence  of conditional expectations to ensure pyramidal
independence.
 
3) Time-homogenity of our processes on the level of continuous
Bernoulli systems can be encoded via the requirement of the existence
of a shift which is compatible with the filtration -- in this case  we
call such a system a non-commutative continuous Bernoulli shift.
These shifts are introduced in \cite{HKK} and provide a
non-commutative extension of Tsirelson's noises or homogeneous
continuous product systems of probability spaces
\cite{Tsir04a}. Similarly, continuous Bernoulli systems are a
non-commutative extension of continuous products   of probability
spaces.
\end{remarks} 

The definition of the notion `continuous Bernoulli system' puts the
whole emphasis on the von Neumann algebras without specifying an
underlying process with independent increments.  In our setting the
information about the increments of the process  will be encoded in
the notion of an additive flow.  As the example of classical Brownian
motion shows, the increments do not need to consist of bounded
operators, thus these flows need not to be elements of the von Neumann
algebras. In general, they will be given by closed densely defined
affiliated operators.   In the following we will restrict to the case
where these flows are   elements of non-commutative $L^p$-spaces, such
that all moments of the   flow exist.
 
\begin{notation} 
Let $\cA$ be a von Neumann algebra and $\ff$ a faithful normal trace.
For $1\leq p<\infty$, the non-commutative $L^p$-spaces $L^p(\cA)$ are
defined as the completion of $\cA$ in the norm
$$\| x\|_p:=\bigl(\ff(\vert x\vert^p)\bigr)^{1/p}\qquad  (x\in\cA),$$
where $\vert x\vert=(x^*x)^{1/2}$.  $L^\infty(\cA)$ is just $\cA$ with
the usual operator norm.  Furthermore, we put
$$L^{\infty-}(\cA):=\bigcap_{1\leq p<\infty}L^p(\cA).$$  Notice that
$\ff$ extends from $\cA$ to $L^1(\cA)$ and that this extension will be
denoted by the same symbol $\ff$.  For further details on
non-commutative $L^p$-spaces we refer to  \cite{PiXu} and the
literature cited therein.
\end{notation} 
 
\begin{definition} 
Let  $(\cA,\ff,(\cA_I)_{I\in\cI})$ be a continuous Bernoulli system.
An \emph{additive flow} (more precisely,   \emph{additive
$L^{\infty-}$-flow}) is a family $B=(B_I)_{I\in\cI}\subset
L^{\infty-}(\cA)$  such that we have
\begin{enumerate} 
\item[(i)]  \emph{continuity}: the map
$$(s,t)\mapsto B_{[s,t)}\in L^p(\cA)$$   is, for all  $1\leq
p<\infty$, jointly continuous  in $s$ and $t$
\item[(ii)]  \emph{adaptedness}: $B_I\in L^{\infty-}(\cA_{I})$ for all
$I\in\cI$
\item[(iii)]  \emph{additivity}: $B_{[s,u)}=B_{[s,t)}+B_{[t,u)}$ for
all $s<t<u$
\end{enumerate} 
 
If $\ff(B_I)=0$ for all $I\in\cI$, then we call the  flow
\emph{centred}. If $B_I\subset\CC 1$  for all $I\in\cI$, then the flow
is \emph{trivial}.  A \emph{normalized} flow is centred and satisfies
$\ff(B_{[0,1)}^*B_{[0,1)})=1$.
 
A flow $(B_I)_{I\in\cI}$ is \emph{stationary} if  we have the
invariance of all its moments in the following  sense: for all
$n\in\NN$ and all $I_1,\dots,I_n\in\cI$ we  have that
$\ff(B_{I_1+t}\dots B_{I_n+t})$ does not depend on $t\in\RR$.
 
For a given stationary flow $(B_I)_{I\in\cI}$ we put
$$B_t:=B_{[0,t)}\qquad (t>0),\qquad B_0:=0$$   and call $(B_t)_{t\geq
0}$ the corresponding  \emph{quantum L\'evy process}.
\end{definition}
 
\begin{remarks} 
1) Notice that we can always turn an additive flow $B_I$ into a
centred additive flow by  considering $B_I-\ff(B_I)$.
 
2) From stationarity and continuity it follows that we have  for a
stationary centered flow that (see also Lemma \ref{lemma:momente})
$$\ff(B_t^2)=\ff(B_1^2)\cdot t.$$  In the same way, by also invoking
the independence of increments,  we get that for any two stationary
centered flows $(B_I)_{I\in\cI}$  and $(\tilde B_I)_{I\in\cI}$ we have
$$\ff(B_{[0,t)}\tilde B_{[0,t)})=\ff(B_{[0,1)}\tilde B_{[0,1)})\cdot
t.$$  In particular, for a normalized stationary flow we have
$$\ff(B_{[0,t)}B_{[0,t)}^*)=t.$$
 
3) Note that we can recover our stationary  flow from the quantum
L\'evy process via
$$B_{[s,t)}=B_t-B_s.$$  This gives $B_I$ only for $I\subset \RR_+$,
however, in the stationary  case this contains all essential
information.  Thus, stationary flows and quantum L\'evy processes are
just two sides of the same object.
\end{remarks} 

In most concrete cases, continuous Bernoulli   systems are given as
von Neumann algebras generated by specified quantum L\'evy
processes. However, there  exist examples of continuous Bernoulli
systems without any non-trivial  quantum L\'evy process (see also
\cite[Theorems 4.4.3 and 6.5.8]{HKK}). In   analogy with the
classification of product systems of Hilbert spaces such   examples
might be addressed as non-type I. We are here mainly interested  in
type I, i.e., those having stationary flows which generate  the von
Neumann algebras. From a probabilistic point of view  it seems to be
appropriate to call such type I continuous Bernoulli  systems
\emph{non-commutative white noises} (see also   \cite[Subsection
6.5]{HKK} for the time-homogeneous setting).  The question of
continuous Bernoulli systems without   non-trivial flows and the
relation of the present  frame   with the work of Arveson \cite{Arv}
and Tsirelson \cite{Tsir04a} on product systems  will be discussed
elsewhere \cite{Kos2}.
 
\section{Isomorphism of continuous Bernoulli systems} 
  
A first canonical problem is  to classify continuous Bernoulli systems
modulo  a notion of isomorphism which respects the filtration.
 
\begin{definition} 
We say that two continuous Bernoulli systems
$(\cA,\ff,(\cA)_{I\in\cI})$ and  $(\cB,\psi,(\cB)_{I\in\cI})$ are
\emph{isomorphic},  if there exists  an isomorphism $\pi:\cA\to\cB$
which respects the filtration, i.e.
$$\pi(\cA_I)=\cB_I\qquad\text{for all $I\in\cI$},$$  and such that
$$\ff=\psi\circ\pi.$$  We will call such a $\pi$ \emph{filtration
preserving}.
\end{definition} 
 
Note that the latter condition on the traces is automatically
fulfilled if the von Neumann algebras $\cA$ and $\cB$ are  factors.
 
This isomorphism problem asks for  a classification of subfactors of
von Neumann algebras in   the extreme case where we have a continuous
family of subalgebras  (of necessarily infinite index).

Since a filtration preserving isomorphism extends to isometries
between the corresponding $L^p$-spaces ($1 \le p < \infty)$,
normalized stationary flows are mapped to normalized  stationary
flows. Thus it is clear that the set of all  distributions of such
flows yields an  invariant for filtration preserving isomorphisms. In
particular,  if we have only one such distribution then this can be
used to  distinguish different continuous Bernoulli systems.

The uniqueness of such a distribution is, for example, given in the
case of classical Brownian motion.  One way to see this is to invoke
the chaos decomposition property of the $L^2$-space of classical
Brownian motion. This says that every element in the $L^2$-space can
be represented (in a unique way) as a sum of multiple Wiener integrals
with respect to Brownian motion.   This means in particular that every
flow can be represented in terms of multiple integrals and by using
the stationarity and the independence of the increments this readily
implies that such a stationary flow has to have Gaussian
distributions.

We want to imitate that argument in the non-commutative case.  The
chaos decomposition of the $L^2$-space into multiple Wiener integrals
equips the $L^2$-space with a Fock space structure $\bigoplus
L^2(\RR^n)$, and the main argument consists then of the simple
observation that non-trivial flows exist in $L^2(\RR^n)$ only  for
$n=1$.
 
As it turns out, in general we do not have such a   chaos
decomposition of the  $L^2$-space of a given continuous Bernoulli
system.   Even if we restrict  to non-commutative versions of Brownian
motions this chaos  decomposition is not present in general. However,
for a quite  big class of continuous Bernoulli systems we have a more
general kind of  chaos decomposition for the corresponding
$L^2$-space,  resembling a Fock space decomposition, but carrying some
additional information.
 
The class of continuous Bernoulli systems for which such a more
general  kind of chaos decomposition is available are the so-called
generalized Brownian motions, which were introduced in \cite{BSp2}.
They are characterized by the requirement that  all mixed moments in
such a Brownian motion can  be calculated by a kind of Wick formula in
terms  of a given function $\ct$ on pair partitions.  In \cite{GM2},
Guta and Maassen have shown that this class  of generalized Brownian
motions coincides with the class  of operators arising in their
construction \cite{GM1} of symmetric  Hilbert spaces in terms of the
combinatorial concept of  species. In particular, they provide a
concrete  realization of the $L^2$-space of the generalized Brownian
motions. Namely, they are of a Fock space like form $\cF_V(\HH)$,
carrying,  however, in general some additional information, which is
encoded  in a sequence $V=(V_n)_{n=0}^\infty$ of (not necessarily
finite  dimensional) Hilbert spaces such that each $V_n$ carries a
unitary  representation $U_n$ of the symmetric group $S(n)$. Then
$$\cF_V(\HH):=\bigoplus_{n=0}^\infty \frac 1{n!} V_n\otimes_s \HH^{
\otimes n},$$   i.e. $\cF_V(\HH)$ is spanned by linear combinations of
vectors of  the form
$$v\otimes_s h_1\otimes\dots\otimes h_n:=\frac 1{n!}  \sum_{\pi\in
S(n)}U_n(\pi)v\otimes \tilde U_n(\pi)h_1\otimes\dots  \otimes h_n,$$
where $\tilde U_n$ is the canonical action of $S(n)$ on the  $n$-fold
tensor product of $\HH$.
 
The concrete structure of the  space $\cF_V$ depends of course on the
scalar product in the spaces  $V_n$, which is determined by the
underlying function $\ct$.  Of course, this Fock space structure is
compatible with the filtration  $I\mapsto L^2(I)$, i.e., under the
identification  of $L^2(\cA)$ with $\cF_V(L^2(\RR))$,   the subspace
$L^2(\cA_I)$ is, for  each $I\in\cI$,  mapped onto $\cF_V(L^2(I))$.
 
It is this form of decomposition for the $L^2$-space  which gives
restrictions for a flow. Although some of the  following arguments
might be extended to more general situations  we will, for sake of
clarity, restrict to the (quite big!) class  of generalized Brownian
motions.  In the following we will denote by $L^2_{\text{loc}}$ the
set of  locally $L^2$-functions, i.e., those measurable functions,
whose  restriction to any compact interval is $L^2$.
 
\begin{theorem} \label{theorem:identification}
Let $(\cA,\ff,(\cA_I)_{I\in\cI})$ be a continuous Bernoulli system,
generated by a generalized Brownian motion, with generalized chaos
decomposition $L^2(\cA)=\cF_V(L^2(\RR))$ for $V=(V_0,V_1,\dots)$.
Then the set of centered flows for $(\cA,\ff,(\cA_I)_{I\in\cI})$   can
be identified with the one-particle space  $V_1\otimes_s
L_{\text{loc}}^2(\RR)\cong L_{\text{loc}}^2(\RR,V_1)$, via
$$L_{\text{loc}}^2(\RR,V_1)\ni\xi\mapsto   (B_I(\xi))_{I\in\cI},$$
where
$$B_I(\xi):=\xi\cdot\chi_I.$$
\end{theorem} 
 
\begin{proof} 
It is clear that any $(B_I(\xi))_{I\in\cI}$ is a centered flow.  (Note
that all moments of these operators exist by the construction  of
generalized Brownian motions and that they are continuous in  the
endpoints of the intervals $I$.)
 
For the other direction,  consider a centered flow $(B_I)_{I\in\cI}$.
Since, by definition, all its moments exist, we must  have that
$B_I\in \cF_V(L^2(I))$.  We decompose  $B_I$ according to the direct
sum decomposition of our $L^2$-space  as
$$B_I=\bigoplus_{n=0}^\infty B_I^{(n)}\qquad\text{with}  \qquad
B_I^{(n)}\in \frac 1{n!} V_n\otimes_s L^2(I)^{\otimes n}  \subset
V_n\otimes L^2(I^n)$$  Note that each $(B_I^{(n)})_{I\in\cI}$ is a
flow, too.  Fix $I\in\cI$ and decompose it, for each $N\in\NN$, into
the  disjoint union of intervals $I_{N,1},\dots,I_{N,N}$ of same
length.  Then, for each $n\in\NN$, we have
$$B_I^{(n)}=B_{I_{N,1}}^{(n)}+\cdots B_{I_{N,N}}^{(n)}\subset
V_n\otimes \bigl(L^2(I_{N,1}^n)\cup\dots\cup L^2(I_{N,N}^n)\bigr).$$
If we send $N\to\infty$, then $B_I^{(n)}$ must live on the
one-dimensional   diagonal in $V_n\otimes L^2(\RR^n)$, which is only
possible for $n=0$  and $n=1$. Centeredness of our flow excludes
$n=0$, and thus we  remain only with the possibility that $B_I\in
L^2(I,V_1)$. Additivity  of the increments yields then that
$B_I=\xi\cdot\chi_I$ for  some locally $L^2$-function $\xi$.
\end{proof} 
 
In many interesting cases, the space $V_1$ is one-dimensional.   In
such a situation a corresponding centered flow must  be of the form
$$B_I=v\otimes f\cdot\chi_I,$$   where $v$ is a fixed unit  vector in
$V_1$ and $f\in L_{\text{loc}}^2(I)$.  If we restrict now to
selfadjoint   stationary normalized flows then we must have
$\ff(B_{[0,t)}B_{[0,t)})=t$ and thus  (note that, because of
selfadjointness, $f$ is real-valued)
$$ t=  \ff(B_{[0,t)}B_{[0,t)})=\la v\otimes
f\cdot\chi_{[0,t)},v\otimes  f\cdot\chi_{[0,t)}\ra=  \int_0^t \vert
f(t)\vert^2 dt,$$  i.e., $f$ must be a function with (almost surely)
constant modulus   $1$.
 
\begin{theorem}  
If the space $V_1$ in Theorem \ref{theorem:identification}  is
one-dimensional then every selfadjoint  stationary normalized flow
$(B_I)_{I\in\cI}$  has the same distribution for  $B_{[0,1)}$.  Thus,
within the class of generalized Brownian motions with one-dimensional
space $V_1$, the distribution for $I=[0,1)$ of the generating flow
$(v\otimes \chi_I)_{I\in\cI}$ is an invariant of the corresponding
continuous Bernoulli systems with respect to filtration preserving
isomorphisms.
\end{theorem} 
 
Note that the distribution of $v\otimes \chi_I$ for arbitray $I$  is
just a dilation of the distribution for $I=[0,1)$, thus does not
contain any additional information.
 
\begin{proof} 
A filtration preserving isomorphism between two  continuous Bernoulli
systems maps a selfadjoint stationary normalized  flow to an object of
the same  kind. For a generalized Brownian motion, the generating flow
$(v\otimes \chi_I)_{I\in\cI}$   is always selfadjoint, stationary and
normalized.   On the other side, as we have seen above,  every
selfadjoint stationary normalized flow must  be of the form $v\otimes
f\chi_I$, where $f$ is a function of  constant modulus $1$. However,
in the calculation of moments for  such operators, only the inner
product between the involved   functions will play a role, which means
that the moments of  $(v\otimes \chi_I)_{I\in\cI}$ are the same as
those of  $(v\otimes f\chi_I)_{I\in\cI}$. Thus the moments of the
generating  flows of two generalized Brownian motions must  be mapped
onto each other by a filtration preserving isomorphism.
\end{proof} 
 
\begin{corollary} 
1) The $q$-Brownian motions (with $-1\leq q\leq 1$) of \cite{BSp1,BKS}
lead  for different $q$ to non-isomorphic continuous Bernoulli
systems.
 
2) The generalized Brownian motions  of \cite{BSp2} lead for different
$q$ to non-isomorphic  continuous Bernoulli systems.
\end{corollary} 
 
\begin{proof} 
Both cases fit into the frame of generalized Brownian motions,  and it
is easy to see that their space $V_1$ is one-dimensional.  Thus  the
distribution of the underlying Brownian motions distinguishes  these
objects with respect to filtration preserving isomorphisms.  It is
easy to see that all distributions are different.
\end{proof}

\section{Moments of quantum L\'evy processes} 
 
Important information about stationary flows $(B_I)_{I\in\cI}$ is
contained in moments of the corresponding quantum L\'evy processes.
 
\begin{lemma} \label{lemma:momente}
Let    $B=(B_I)_{I\in\cI}$ be a stationary flow and $(B_t)_{t\geq 0}$
the corresponding quantum L\'evy process.   Then there exist constants
$\alpha$, $\beta$, and $\gamma$ such that  we have for all $t>0$
\begin{align*} 
\ff(B_t)&=\alpha t,\\  \ff(B_t^2)&=\alpha^2 t^2 +\beta t,\\
\ff(B_t^3)&=\alpha^3 t^3+ 3\alpha\beta t^2+ \gamma t.
\end{align*} 
\end{lemma} 
 
\begin{proof} 
For all $s,t\geq 0$, we have
$$B_{s+t}=B_{[0,s)}+B_{[s,s+t)},$$  and thus
$$\ff(B_{s+t})=\ff(B_s)+\ff(B_t),$$  which gives, by continuity, the
equation for the first  moment, with $\alpha=\ff(B_1)$.
 
For the second moment we get
$$B_{s+t}^2=B_{[0,s)}^2+B_{[0,s)}\cdot B_{[s,s+t)}+B_{[s,s+t)}\cdot
B_{[0,s)}+B_{[s,s+t)}^2.$$  Pyramidal independence gives
$$\ff(B_{s+t}^2)=\ff(B_s^2)+\ff(B_t^2)+2\ff(B_s)\ff(B_t),$$  which
implies the equation for the second moment.
 
Similarly, one shows the result for the third moment.
\end{proof} 
 
Note that pyramidal independence does not allow to  calculate all
mixed moments of fourth and higher order: e.g.,  we cannot make a
general statement about
$\ff(B_{[0,s)}B_{[s,s+t)}B_{[0,s)}B_{[s,s+t)})$.   Thus, in this
generality, similar statements   as in the Lemma \ref{lemma:momente}
are not accessible for higher moments.  Nevertheless, explicit
polynomial bounds for the growth of   higher moments are established
in \cite{Kos,Kos4,Kos1}, as an application of   Burkholder-Gundy
resp. Burkholder/Rosenthal inequalities for non-commutative
$L^p$-martingales \cite{PX,JX}.
 
However, if we require some more special structure, then we can say
much more  about the behaviour of higher moments.   In this section we
want to consider the case where we have  an order invariance of the
moments of the increments, in the sense  that such moments do not
change if we shift the increments against  each other as long as we do
not change the relative position of  the intervals.   Let us first
consider a discrete version of this before we treat  the continuous
case.
 
\subsection{Limit theorem for order invariant distributions} 
 
Consider random variables $b_i^{(N)}$ ($i,N\in\NN$, $i\leq N$) living
in some non-commutative probability space $(\cA,\ff)$.
 
For an $n$-tuple
$$\ii:\{1,\dots,n\}\to\{1,\dots,N\}$$  we put
$$b_\ii^{(N)}=b_{i(1)}^\hN\cdots b_{i(n)}^\hN.$$  For an $\ii$ as
above, we denote by $\vert \ii\vert$ the  number of elements in the
range of $\ii$.
 
\begin{definition} 
1) Let $\ii,\jj:\{1,\dots,n\}\to \NN$ be two $n$-tuples of
indices. We say that they are \emph{order equivalent}, denoted by
$\ii\sim \jj$, if
$$i(k)\leq i(l)\Longleftrightarrow j(k)\leq j(l)  \qquad\text{for all
$k,l=1,\dots,n$.}$$  We denote by $\cV(n)$ the set of equivalence
classes for maps  $\ii:\{1,\dots,n\}\to\NN$ under this order
equivalence. Note  that for each $n$ this is a finite set.
 
2) We say that the distribution of the variables $b_i^\hN$ is
\emph{order invariant} if we have for all $n,N\in\NN$ and all
$\ii,\jj:\{1,\dots,n\}\to\{1,\dots,N\}$ with $\ii\sim \jj$  that
$$\ff(b_\ii^\hN)=\ff(b_{\jj}^\hN).$$  In this case we denote, for
$\sigma\in\cV(n)$,   by $\ff(b_\sigma^\hN)$ the common value  of
$\ff(b_\ii^\hN)$ for $\ii\in\sigma$.
\end{definition} 
 
Given such order invariant random variables, one can make quite
precise statements about the moments of the sums
$b_1^\hN+\dots+b_N^\hN$ in the limit $N\to\infty$. The proof  of this
limit theorem follows the usual arguments, see, e.g.,  \cite{SpW}, and
we will omit the proof.
 
\begin{theorem} \label{limit-thm} 
Consider random variables $b_i^{(N)}\in(\cA,\ff)$ ($i,N\in\NN$, $i\leq
N$),   whose distribution is order invariant.  Assume that for all
$n\in\NN$ and all  $\sigma\in\cV(n)$ the following limit exists:
$$c(\sigma):=\lim_{N\to\infty} N^{\vert
\sigma\vert}\ff(b_{\sigma}^\hN)  .$$  Define
$$S_N:=b_1^\hN+\dots+b_N^\hN.$$  Then we have for all $n\in\NN$
$$\lim_{N\to\infty}\ff(S_N^n)=\sum_{\sigma\in\cV(n)}\alpha_\sigma
c(\sigma),$$  where the $\alpha_\sigma$ are the constants,
$$\alpha_\sigma=\lim_{N\to\infty}\frac{\#\{\ii:\{1,\dots,n\}\to
\{1,\dots,N\}\mid \ii\in\sigma\}}{N^{\vert\sigma\vert}}
=\frac1{\vert\sigma\vert!}.$$
\end{theorem} 
 
\subsection{Moments of order invariant quantum L\'evy processes} 
In the following, we will use,  for two intervals $I,J\in\cI$, the
notation $I<J$   to indicate that we have $s<t$ for all $s\in I$ and
$t\in J$.
 
\begin{definition} 
Let $(B_I)_{I\in\cI}$ be a flow. We say that the flow (or its
corresponding quantum L\'evy process)  is \emph{order invariant} if
we have for all $I_1,\dots,I_n\in\cI$ with $I_k\cap I_l=\emptyset$
($k,l=1,\dots,n$) that
$$\ff(B_{I_1}\cdots B_{I_n})=\ff(B_{I_1+t_1}\cdots B_{I_n+t_n})$$  for
all $t_1,\dots,t_n$ with the property that, for all $k,l=1,\dots,n$,
$I_k<I_l$ implies $I_k+t_k<I_l+t_l$.
\end{definition} 
 
\begin{remark} 
Note that an order invariant flow is in particular stationary.
\end{remark} 
 
Consider now such an order invariant   flow $(B_I)_{I\in\cI}$.  Put
$$b_i^\hN:=B_{[\frac{i-1}N,\frac iN)}.$$  Then we have
$$S_N=b_1^\hN+\cdots + b_N^\hN=B_1$$  for all $N\in\NN$ and, since the
distribution of the $b_i^\hN$ is order invariant, our Limit  Theorem
\ref{limit-thm} yields that
$$\ff(B_1^n)=\sum_{\sigma\in\cV(n)}\alpha_\sigma c(\sigma),$$  if all
$$c(\sigma):=\lim_{N\to\infty} N^{\vert
\sigma\vert}\ff(b_{\sigma}^\hN)$$  exist.
 
\begin{proposition} 
Let $(B_I)_{I\in\cI}$ be an order invariant flow. Then, for  all
$n\in\NN$ and $\sigma\in\cV(n)$, the limit
$$c(\sigma):=\lim_{N\to\infty} N^{\vert
\sigma\vert}\ff(b_{\sigma}^\hN)$$  exists.
\end{proposition} 
 
\begin{proof} 
 
We will prove this, for fixed $n$,   by induction over the length of
$\sigma$, starting  with maximal length of $\sigma$.
 
Namely, fix $n$ and consider first  a $\sigma$ with
$\vert\sigma\vert=n$. This  means that
$\ii=(i(1),\dots,i(n))\in\sigma$ is a tuple of $n$  different
numbers. By using the stochastic independence we get
\begin{align*} 
N^{\vert \sigma\vert}\ff(b_{\sigma}^\hN)&=  N^n\ff(b_{i(1)}^\hN\cdots
b_{i(n)}^\hN)\\  &= N^n\ff(b_{i(1)}^\hN)\cdots\ff(b_{i(n)}^\hN)\\
&=N^n\ff(b_1^\hN)^n  \\&=\bigl(N\ff(B_{[0,\frac 1N)})\bigr)^n  \\&=
\ff(B_1)^n,
\end{align*} 
and hence the limit
$$c(\sigma):=\lim_{N\to\infty} N^{\vert
\sigma\vert}\ff(b_{\sigma}^\hN)  =\ff(B_1)^n$$  exists.
 
Consider now an arbitrary $\sigma\in\cV(n)$ and assume that we have
proved the existence of the limits $c(\sigma')$ for all $\sigma'\in
\cV(n)$ with $\vert\sigma'\vert>\vert\sigma\vert$.  Choose an
$n$-tuple $\ii=(i(1),\dots,i(n))\in\sigma$ and consider
$$\ff(B_{[i(1),i(1)+1)}\cdots B_{[i(n),i(n)+1)}).$$  Again, we
decompose the intervals of length 1 into $N$ subintervals  of length
$1/N$, so that we can write this also as
$$\ff\bigl((\sum_{k(1)=1}^NB_{[i(1)+\frac {k(1)-1}N,i(1)+\frac {k(1)}
N)})  \cdots   (\sum_{k(n)=1}^NB_{[i(n)+\frac {k(n)-1}N,i(n)+\frac
{k(n)}N)})  \bigr).$$  If we multiply this out and collect terms
together with the same  relative position of the subintervals then we
get a sum of terms,  one of which is exactly
$N^{\vert\sigma\vert}\ff(b_\sigma^\hN)$, and the others are of the
form  $\gamma_{\sigma'}\ff(b_{\sigma'}^\hN)$, for $\sigma'$ with
$\vert\sigma'\vert>\vert\sigma\vert$. Since also $\gamma_{\sigma'}\sim
N^{\vert\sigma'\vert}$ for $N\to\infty$, we know by our induction
hypothesis that all these other terms have a finite limit for
$N\to\infty$. Since the left hand side of our equation does not
depend on $N$, also the term   $N^{\vert\sigma\vert}\ff(b_\sigma^\hN)$
must have a finite limit  for $N\to\infty$.
\end{proof} 
 
Of course, the same argument works if we replace the time $1$ by  an
arbitrary time $t$. In this case, we get the existence of the  limits
$$c_t(\sigma):=\lim_{N\to\infty}N^{\vert\sigma\vert}
\ff\bigl(B_{[i(1),i(1)+t/N)}\cdots B_{[i(n),i(n)+t/N)}\bigr),$$  for
$(i(1),\dots,i(n))\in\sigma$. The remaining question is  how these
$c_t(\sigma)$ depend on the time $t$.
 
\begin{lemma} 
We have that
$$c_s(\sigma)=c(\sigma)\cdot s^{\vert\sigma\vert}
 \qquad\text{for
all $s\in\QQ$}.$$
\end{lemma} 
 
\begin{proof} 
For $(i(1),\dots,i(n))\in\sigma$ and $t\in\RR$, we have
\begin{align*} 
c_{2t}(\sigma)&=\lim_{N\to\infty}N^{\vert\sigma\vert}
\ff\bigl(B_{[i(1),i(1)+2t/N)}\cdots B_{[i(n),i(n)+2t/N)}\bigr)\\
&=\lim_{N\to\infty}N^{\vert\sigma\vert}
\ff\bigl((B_{[i(1),i(1)+t/N)}+  B_{[i(1)+t/N,i(1)+2t/N)})\cdots
\\&\qquad\qquad\qquad\cdots (B_{[i(n),i(n)+t/N)}
+B_{[i(n)+t/N,i(n)+2t/N)})\bigr)\\ 
&=2^{\vert\sigma\vert}c_t(\sigma).
\end{align*} 
Note that for each block of $\sigma$ we can choose either
 the
increments from $i$ to $i+t/N$ or the increments from
 $i+t/n$ to
$i+2t/N$ to match up, i.e., each block of
 $\sigma$ contributes a
factor $2$. On the other hand, terms which
 match for some block 
an increment from $i$ to $i+t/N$ with an increment from
 $i+t/N$ to
$i+2t/N$  vanish in  the limit, because they correspond to a $\sigma'$
with   $\vert\sigma'\vert>\vert\sigma\vert$, and so they have to be
multiplied with a higher power of $N$ to give a non-trivial  limit.

 In the same way as above one can also see that for any $k\in\NN$
and any $t\in\RR$ we have
$$c_{kt}(\sigma)=k^{\vert\sigma\vert}c_t(\sigma).$$
 This yields
finally the assertion.
\end{proof} 
 
 By invoking different $t$ for each block of $\sigma$ one could
also derive functional equations for these quantities which,
together with the fact that they are measurable, would extend the
statement of the above lemma to all $t\in\RR$. However, we do not
need this because the continuity of the moments $\ff(B_t^n)$ allows
us to extend the statement in the next theorem directly from
rational
 to all real times $t$.
 
 Let us summarize in the
following theorem our results.
 
\begin{theorem} 
Let $(B_t)_{t\geq 0}$ be an order invariant quantum L\'evy process.
Then there exist constants $c(\sigma)$ for all $\sigma\in\cV$  such
that we have for all $n\in\NN$ and all $t\geq 0$

$$\ff(B_t^n)=\sum_{\sigma\in\cV(n)}\frac 1{\vert\sigma\vert!}
c(\sigma) t^{\vert  \sigma\vert}.$$
\end{theorem} 
 
In the next section we will see that quantum L\'evy processes which
are not  order invariant do not necessarily have such a polynomial
behaviour of their moments.
 
\subsection{An example of a non order invariant generalized Brownian
motion}
\label{subsection:example-non-order} 
Finally, we want to present an example   of a quantum L\'evy process
which is not order invariant.  This example is a generalization of the
$q_{ij}$-relations
$$a_ia^*_j-q_{ij}a^*_ja_i=\delta_{ij}\cdot 1$$   to the continuous
case, and it is formally given by
$$a_ta^*_s-q(s-t)a^*_sa_t=\delta(s-t)\cdot 1.$$
 
This situation can be realized rigorously as follows:  Put
$\HH:=L^2(\RR)$, and consider on $\HH\otimes\HH=  L^2(\RR^2)$ the
operator $T$, given by
$$(Tf)(s,t)=q(s,t)\cdot f(t,s),$$  where $q=q(\cdot,\cdot)$ is a fixed
function of two variables.  This $T$ fulfills the braid relations. If
we assume  in addition that $q$ has the properties $\bar
q(s,t)=q(t,s)$ and
 $\vert q(s,t)\vert\leq1$ for all $s,t$,   then
$T$ is also   selfadjoint and contractive. Thus the  assumptions of
\cite{BSp3} are fulfilled and the corresponding  Fock space
construction yields a positive inner product   and, for each
$f\in\HH$, creation and annihilation operators $d^*(f)$ and  $d(f)$,
respectively.   Put now, for $I\in\cI$,
$$B_I:=d(\chi_I)+d^*(\chi_I)$$  and define
\begin{align*} 
\cA:&=\vN(B_I\mid I\in\cI)\\  \cA_I:&=\vN(B_J\mid J\in\cI, J\subset
I)\qquad (I\in\cI)\\  \ff(a):&=\la a \Omega,\Omega\ra\qquad
\qquad\qquad\qquad(a\in\cA)
\end{align*} 
If $q$ is real (and thus symmetric, i.e., $q(s,t)=q(t,s)$), then $\ff$
is  a faithful trace on $\cA$.  Furthermore, if $q$ is stationary,
i.e. $q(s,t)=q(s-t)$, then one  has a well-defined second quantization
$\Gamma(S_t)$   (see \cite{Kro}) of the  usual shift ($u\in\RR$)
$$S_u: L^2(\RR)\to L^2(\RR),\qquad (S_u f)(t)=f(t-u),$$  which is
compatible with the filtration of the von Neumann algebras.  Let us
summarize this in the following proposition.
 
\begin{proposition} 
Let $q:\RR\to\RR$ be a measurable function with the property
$$-1\leq q(t)=q(-t)\leq 1\qquad \text{for all $t\in\RR$},$$  then
$(\cA,\ff,(\cA_I)_{I\in\cI}$  corresponding to the operator $T$ on
$L^2(\RR^2)$ given by
$$(Tf)(s,t)=q(s-t)\cdot f(t,s)$$  forms a continuous Bernoulli system
with corresponding stationary  flow
$$B_I:=d(\chi_I)+d^*(\chi_I) \qquad (I\in\cI).$$
\end{proposition} 
 
If $q$ is constant, then   one recovers the example of the
$q$-Brownian motion \cite{BSp1,BKS},  which is of course order
invariant.  If, however, $q$ is not constant then this flow is not
order invariant.   For example, by using the definition of the
operators $d(f)$  and $d^*(f)$, one readily finds for $I,J\in\cI$ with
$I\cap J  =\emptyset$ that
$$\ff(B_IB_JB_IB_J)=\int_I\int_J q(s-t)dsdt,$$  which gives for the
fourth moment of our quantum L\'evy process
$$ 
\ff(B_t^4)=t^2+\int_0^t q(t_1-t_2)dt_1dt_2.
$$ 
 
(Note that formally these results can be obtained by using the  Ito
formula
$$dB_s dB_t dB_s dB_t =q(s-t)ds dt.)$$
 
By making different choices of the function $q$, this shows that
there is quite a variation of the behaviour  of the fourth (and
higher) moments for non order invariant quantum L\'evy  processes.

\end{document}